\documentstyle{amsppt}
\input amstex
\magnification 1100
\hsize=15 truecm
\vsize=21.5 truecm
\NoBlackBoxes

\topmatter
\title COMBINATORIAL METHODS:  FROM GROUPS TO POLYNOMIAL ALGEBRAS 
\endtitle
\author  Vladimir Shpilrain 
\endauthor 

\abstract{\baselineskip 12pt
  Combinatorial methods (or methods of elementary
transformations)   came to group theory from low-dimensional
topology in the beginning of  the century. Soon after that,
combinatorial group theory   became an  independent area with its
own powerful techniques. On the other hand,  combinatorial
commutative algebra emerged in the sixties, after Buchberger 
introduced what is now known as Gr\"{o}bner bases.   The purpose of
this survey  is to show how ideas from one of those areas contribute
to the other. }

\endabstract

\endtopmatter
\document
\baselineskip 15pt 

\heading 1. Introduction\endheading

   Let  $F = F_n$   be the free group of a finite rank  
$n \ge 2$  with a  set  
$X = {\{}x_1,...,x_n {\}}$  of free generators.
  Let   $Y  =  {\{}y_1,...,y_m{\}}$   
 and  $\widetilde Y = {\{}\widetilde y_1,...,\widetilde y_m{\}}$  
be arbitrary finite sets of elements of the  group  $F$.  Consider 
the  following  elementary  transformations  that  can  be applied
to  $Y$:  
\smallskip 

\noindent {\bf (N1)} $y_i$   is replaced by  $y_i y_j$ or by  
$y_j y_i$ for some  $j \ne i$ ;

\noindent {\bf (N2)} $y_i$   is replaced by  $y_i^{-1}$. 

\noindent {\bf (N3)} $y_i$   is replaced by some $y_j$,  and  at
the  same time $y_j$   is replaced by $y_i$. 

\smallskip 

    It is understood that  $y_j$  doesn't change if  $j\ne i$. 
\smallskip 

 One might notice that some of these transformations  are 
redundant, i.e., are compositions of other ones. There is 
a reason behind that which we are going to explain a little later. 
\smallskip 

    We say that two sets   $Y$  and   $\widetilde Y$  are Nielsen
equivalent  if one of them can be obtained from another by  applying 
a sequence of transformations (N1)--(N3). It was proved by Nielsen 
 that two sets 
$Y$  and   $\widetilde Y$ generate the same subgroup of the 
group  $F$    if and  only if they are Nielsen equivalent. This 
result is now one of the central points in combinatorial group
theory.  
\smallskip 

 Note however that this result alone does not give an {\it
algorithm} for deciding whether or not $Y$  and   $\widetilde Y$ 
generate  the same subgroup of $F$. To obtain an algorithm, 
we need to somehow define the {\it complexity} of a given set of 
elements, and then to show that a sequence of Nielsen 
 transformations (N1)--(N3) can be arranged so that this complexity 
 decreases (or, at least, does not increase) {\it at every step}
(this is where we may need ``redundant" elementary  transformations
!).

 This was also done by Nielsen; the complexity of a given set 
$Y  =  {\{}y_1,...,y_m{\}}$   is just the sum of the 
 lengths of the
words  $y_1,...,y_m$. We refer to [15] for details. 

 Nielsen's method therefore yields (in particular) an algorithm 
for deciding whether or not a given endomorphism  of a free group
of finite rank is
actually an automorphism. 
\smallskip 

  A somewhat more difficult problem is, given a pair of 
 elements of a free group $F$, to find out if one of them can be 
taken to another by an automorphism of $F$. We call this problem 
 {\it the automorphic conjugacy problem}. It was addressed by 
Whitehead who came up with another kind of elementary  
transformations in a free group: 
\smallskip 

\noindent {\bf (W1)} For some $j$, every 
$x_i$,  $i \ne j$,   is replaced by  one of the elements 
$x_i x_j$,    $x_j^{-1} x_i$, $x_j^{-1} x_i x_j$, or $x_i$. 

\noindent {\bf (W2)} $x_i$   is replaced by  $x_i^{-1}$. 

\noindent {\bf (W3)} $x_i$   is replaced by some $x_j$,  and  at
the  same time $x_j$   is replaced by $x_i$. 
\smallskip 

 One might notice a similarity of Nielsen and Whitehead 
 transformations. However, they differ in one essential detail: 
Nielsen transformations are applied to arbitrary sets of elements, 
whereas Whitehead transformations are applied to a {\it fixed
basis} of the group $F$. 
\smallskip 

Using (informally)  matrix language, 
we can say that Nielsen transformations correspond to elementary  
rows 
 transformations of a matrix (this correspondence can actually 
be made quite formal -- see [17]), whereas Whitehead transformations 
 correspond to  conjugations (via changing the  basis). This latter 
type of matrix transformation is known to be more complex, and 
the corresponding structural results are deeper. 

 There is very much the same relation   between Nielsen and
Whitehead  transformations in a free group. 
\smallskip 

 Note also that Whitehead transformation (W1) is somewhat more 
complex than its analog (N1). This is -- again -- in order 
 to be able to 
arrange a sequence of elementary transformations so that 
 the   complexity of a given element (in this case, just the 
lexicografic length of a cyclically reduced word) would 
decrease (or, at least, not increase) {\it at every step} -- see
    [14]. 

 This arrangement still leaves us with a difficult problem - 
to find out if one of two  elements {\it of the same 
   complexity} (= of the same length) can be 
taken to another by an automorphism of $F$. This is actually 
the most difficult part of Whitehead's algorithm. 
\smallskip 

 In one special case however this problem does not arise, namely, 
when one of the elements is {\it primitive}, i.e., is an
automorphic image of $x_1$. If we have managed to reduce an element
of  a free group (by Whitehead transformations) to an element of 
 length 1, we immediately conclude that it is primitive; no further 
 analysis is needed. 

 Thus, the problem of distinguishing primitive elements of 
a free group is a relatively easy case of the automorphic conjugacy 
problem. As we shall see in Section 2, this is also  the 
situation  in a polynomial algebra. 
\medskip

 In Section 2, we review the results of various attempts to 
create something similar to Nielsen's and Whitehead's methods 
for a  polynomial algebra in two variables. For a  polynomial 
algebra in more than two variables, these problems are
unapproachable so far, since we don't even know what the generators 
of the automorphism group of an algebra like that  look like.
\smallskip 

  In Section 3, we  talk about {\it retracts} of a  polynomial 
algebra in two variables. Basic properties of retracts of a free 
group are given in [15].  Since then, retracts have not been 
getting much attention until very recently, when Turner [21] and, 
independently, Bergman [4]   brought them back to life by 
 employing  them in various interesting research projects 
in combinatorial group theory. Here we show the relevance of 
polynomial retracts to several well-known problems about 
 polynomial mappings, in particular, to the notorious  Jacobian 
conjecture.  
\smallskip 

 In the concluding Section 4, we have gathered some open 
combinatorial problems about polynomial mappings that are motivated
by similar issues in combinatorial group theory. 
\smallskip

\heading 2. Elementary  transformations  in polynomial 
algebras  \endheading 

  Let $P_n = K[x_1,...,x_n]$ be the polynomial algebra in   $~n$
 variables over a field $K$ of characteristic $0$.
 We are going to concentrate here mainly on the algebra $P_2$. 

 The first description of the group $Aut(P_2)$ was given by Jung
[12] back 
 in 1942, but it was limited to the case $K = {\Bbb
C}$ 
 since  he was using methods of algebraic geometry. Later on, van
 der  Kulk   extended  Jung's result to arbitrary ground fields.
 In the form we give it here, the result  appears as Theorem 8.5 in
P.M.Cohn's book [6]; this form is  consistent with the idea of
elementary  transformations as described in the Introduction. 

\proclaim{Theorem 2.1}[6] Every automorphism of $K[x_1,x_2]$ is 
a product of linear automorphisms and automorphisms of the form 
$x_1 \to x_1 + f(x_2); ~x_2 \to x_2$. More precisely, if $(g_1,
g_2)$ is an automorphism of   $K[x_1,x_2]$  such 
  that  $~deg(g_1) \ge deg(g_2)$, say, then either $(g_1, g_2)$ is a
linear automorphism, or there exists a unique $~\mu \in 
 K^{\ast}$  and a positive integer $~d~$  such 
  that  $~deg(g_1 - \mu g_2^d) < deg(g_1)$. 
\endproclaim 

The proof given in [6] is attributed to Makar-Limanov
(unpublished), 
 with simplifications by Dicks [8]. 
\smallskip 

 Note that the ``More precisely, ..."  statement serves the 
algorithmic purposes: upon defining the complexity of a given 
pair of polynomials $(g_1, g_2)$ as the sum $~deg(g_1) + deg(g_2)$, 
we see that Theorem 2.1 allows one to arrange a sequence of 
 elementary  transformations (these are linear automorphisms and 
automorphisms of the form  $x_1 \to x_1 + f(x_2); ~x_2 \to x_2$) 
 ~so that this complexity decreases  
 at every step, until we either get a pair of polynomials that 
represents a linear automorphism, or conclude that $(g_1, g_2)$ 
was not an automorphism of   $K[x_1,x_2]$. The parallel  with 
 Nielsen's method described in the Introduction is obvious. 
\smallskip 

 We also mention here another proof of this result (in case $char
~K = 0$) due to 
Abhyankar and Moh  [1]. In fact, their method is even more 
similar to  Nielsen's method in a free group. Many of their results
are  based on the following fundamental theorem 
which we give here only in the characteristic $0$ case 
(it will also play an essential role in our Section 3): 

\proclaim{Theorem 2.2}[1] Let $~u(t), v(t) \in K[t]$ be two 
 one-variable polynomials of degree $n \ge 1$ and $m \ge 1$. 
 Suppose $~K[t] = K[u,v]$. Then either $~n$ divides $~m$, or 
 $~m$ divides $~n$. 
 \endproclaim 
\smallskip 

 Now let's see how one can adopt a more sophisticated Whitehead's 
method in a polynomial algebra situation. It appears that 
 elementary basis  transformations (see Theorem 2.1), when applied 
 to a polynomial $p(x_1, x_2)$, are mimicked by 
 Gr\"{o}bner transformations of a basis of the ideal of $P_2$
generated by 
 partial derivatives of this polynomial. To be more specific, we 
have to give some background material first. 
\smallskip 

In the course of constructing a Gr\"{o}bner basis 
of a given ideal of $P_n$, one uses ``reductions'', i.e.,
transformations of the following type (see [2], p.39-43): given a
pair $(p, q)$ of  polynomials, set ~$ S(p, q) = {L \over l.t.(p)}~p -
{L \over l.t.(q)}~q$, ~~where ~$l.t.(p)$ ~is the {\it leading
term} of ~$p$, i.e., the {\it leading monomial} together with its
coefficient; ~$L = l.c.m.(l.m.(p), l.m.(q))$ (here, as usual,
$l.c.m.$ means the least common multiple, and $l.m.(p)$ denotes the
 leading  monomial of ~$p$). In this paper, we'll always consider
what is called ``deglex  ordering'' in [2] - where monomials are 
ordered first by  total degree, then lexicographically with 
$x_1 > x_2 > ... > x_n$. 
\medskip

Now a crucial observation is as follows. These Gr\"{o}bner
reductions appear to be of two essentially different types: 
\smallskip 

(i) {\it regular}, or {\it elementary}, transformations. These are
of the form $ S(p, q) = \alpha \cdot p - r \cdot q$ $~$~or $~S(p, q)
= \alpha \cdot q  - r \cdot p$ ~ ~for some polynomial ~$~r$ and
scalar $\alpha \in K^*$. This happens when the
 leading  monomial of ~$~p$ ~is divisible by the leading  monomial of
~$q$ ~(or vice versa).  The reason why we call these transformations 
{\it elementary} is that they can be written in the form $(p, ~q)
\to  (\alpha_1 p, ~\alpha_2 q) \cdot M$, ~where $M$ is an {\it
elementary matrix}, i.e., a matrix which (possibly) differs from the
identity matrix by a single element outside the diagonal. In case
when we have more than 2  polynomials ~$(p_1, ..., p_k)$, we also
can write  $(p_1, ..., p_k) \to  (\alpha_1 p_1, ..., \alpha_k
p_k) \cdot M$, $~$~where $M$ is a $k \times k$ elementary matrix;  
elementary  reduction here is actually applied to  a pair of
polynomials  (as usual)  while the other ones are kept fixed.
Sometimes, it is more convenient for us to get rid of the
coefficients $\alpha_i$  and   write $(p_1, ..., p_k) \to   ( p_1,
..., p_k) \cdot M$, ~where $M$ belongs to the group $GE_k(P_n)$
generated by  all elementary {\bf and } diagonal matrices from
$GL_k(P_n)$. It is known [20] that $GE_k(P_n) = GL_k(P_n)$ if $k \ge
3$, and  $GE_2(P_n) \ne GL_2(P_n)$ if $n \ge 2$ - see [5]. 
 \smallskip 

(ii) {\it singular } transformations -- these are non-regular ones. 
\medskip

  Denote by $~I_{d(p)}$ the 
ideal   of ~$P_2$ ~generated by partial derivatives of $~p$. 
We say that a polynomial $~p \in P_n$ has a {\it unimodular gradient}
if $~I_{d(p)} = P_n$ (in particular, the ideal $~I_{d(p)}$ has rank
1  in this case). Note that if the ground field  $K$ is 
algebraically closed, then this  is equivalent, by Hilbert's
Nullstellensatz, to the  gradient being nowhere-vanishing. 
\smallskip 
 
 Furthermore, define the {\it outer rank} of a polynomial  $p \in
P_n$ to be the  minimal  number   of generators ~$x_i$ ~on which an
automorphic  image  of  $ p$   can  depend. 

Then we have:

\proclaim{Theorem 2.3} [18] Let a polynomial $~p \in P_2$ have
unimodular  gradient.  Then the outer
rank of $~p$ equals 1 if and only if one can get from $(d_1(p),
~d_2(p))$ to  ~(1, ~0)  ~by using only elementary
transformations. ~Or, in the matrix form: if and only if 
~$~(d_1(p), ~d_2(p)) \cdot M = (1, ~0)$ for some matrix $~M  \in 
GE_2(P_2)$. \endproclaim 

  The proof [18] of Theorem 2.3 is based on a generalization of
Wright's Weak Jacobian Theorem [22]. 
\medskip 

\noindent {\bf Remark 2.4.} Elementary transformations that reduce 
$(d_1(p), ~d_2(p))$ to  ~(1, ~0), $~$ can be actually chosen to be 
Gr\"{o}bner reductions, i.e., to decrease the maximum degree of 
monomials {\it at  every step} -- the  proof [18] is based on a
recent result of Park [16]. 
\smallskip 

 Now we show how one can apply this result to the study of 
so-called coordinate polynomials. 

We call a polynomial ~$p \in P_n$ {\it ~coordinate} ~if it can  be
included in  a generating set
of cardinality $~n$ of the algebra $P_n$. It is clear that 
the outer rank of a coordinate polynomial equals 1 (the converse 
 is not true!).  It is easy to show  that a coordinate polynomial
has  a unimodular gradient, and again -- the converse  is not true!
On the other hand, we have: 

\proclaim{Proposition 2.5}[18] A polynomial $p \in P_n$ is
 coordinate  if and only if it has outer rank 1  $~${\bf  ~and } a
unimodular gradient.  
\endproclaim 

Combining this proposition with Theorem 2.3 
 yields the following

\proclaim{Theorem  2.6}[18] A polynomial $p \in P_2$ is coordinate 
if and only if one can get from $(d_1(p), ~d_2(p))$ to  ~(1, ~0)  
~by using only elementary Gr\"{o}bner reductions. 
\endproclaim 

 This immediately yields an algorithm for detecting 
coordinate polynomials in $P_2$ (see [18]) 
 which is similar to  Whitehead's 
 algorithm for detecting primitive elements in a free group. 
This  algorithm is very simple and fast: it has quadratic  growth
with respect  to the degree of a polynomial.      In case $~p$ is
revealed to be a coordinate polynomial, the algorithm also gives a 
polynomial which completes  $~p$ to a basis of $P_2$. 
\smallskip 

In the case when $K = {\Bbb C}$,  the
field of complex numbers,  an 
 alternative, somewhat
more  complicated algorithm, has been recently reported 
 in [9]. It is
not known  whether or not there is an algorithm for detecting 
coordinate polynomials in ~$~P_n$ ~if $~n \ge 3$. 
\medskip 

 Theorems 2.3 and  2.6 also suggest the following conjecture which
is  relevant to an  important problem known as ``effective  
  Hilbert's Nullstellensatz": 
\medskip 

\noindent {\bf Conjecture ``G".}  Let a polynomial $~p \in P_2$ 
have a unimodular  gradient. Then  one can get from $(d_1(p),
~d_2(p))$ to  ~(1, ~0)  ~by using {\it at most one}  singular 
Gr\"{o}bner  reduction. 
 
\medskip 

\noindent {\bf Remark 2.7.} For $n \ge 3$, $~$ Theorem 2.3  is no
longer valid since in this case, by a result of Suslin [20], the 
group $~GL_n(P_n) = GE_n(P_n)$ acts  transitively on  the  set  of 
all  unimodular polynomial vectors of dimension $n$, yet there are 
polynomials with unimodular gradient, but of the outer rank 2, for 
example, $~p = x_1 + x_1^2 x_2$. The ``only if'' part however is
valid for an arbitrary $n \ge 2$ - see  [18].  It is also easy to
show  that  one always has $~orank~p \ge rank(I_{d(p)})$.  
\medskip 

Finally, we mention 
 that our method also yields   an algorithm which, 
given a coordinate  polynomial $~p \in P_2$, finds a sequence  of
elementary  automorphisms (i.e.,  automorphisms of the form $x_1\to
x_1 + f(x_2); ~x_2 \to x_2~ $ together with linear automorphisms)
that  reduces $~p~$ to $~x_1~$. 
\smallskip

\heading 3. Polynomial retracts \endheading

 Let $ K[x, y]$ be the polynomial algebra in two variables 
 over a field $K$ of characteristic $0$. A subalgebra $R$ of  
$K[x, y]$ is called a  {\it retract} if it satisfies any of the
following equivalent conditions:  
\medskip

 \noindent {\bf (R1)} There is an idempotent homomorphism  (a {\it
retraction}, or {\it projection}) 
 $\varphi:   K[x, y] \to  K[x, y]$    such 
  that   $\varphi(K[x, y]) = R.$ 
\smallskip

 \noindent {\bf (R2)} There is a homomorphism   $\varphi:         
K[x, y] \to R$ that fixes every element of  $R$. 
\smallskip 

 \noindent {\bf (R3)}  $ K[x, y] = R \oplus I$  ~for some  ideal   
$I$ of the algebra   $K[x, y]$. 
\smallskip 

 \noindent {\bf (R4)} $ K[x, y]$ is a projective extension of   $R$
in the category of  $ K$-algebras. 
   In other words, there is a splitting exact sequence   $1  \to  I
\to  K[x, y] \to R \to 1$, 
 where  $I$ is the same ideal as in (R3) above. 

\medskip

 Examples: $ K$; $ K[x, y]$; any  subalgebra of the form $K[p]$, 
where $p \in  K[x, y]$ is a  {\it coordinate}  polynomial (i.e., 
$K[p, q] = K[x, y]$ for some  polynomial  $q \in  K[x, y]$). 
 There are other, less obvious, examples of retracts: if  $~p = x +
x^2y$,  then  $K[p]$ is  a retract of  $ K[x, y]$, but $p$ is not 
coordinate since it has a fiber 
  ${\{}p = 0{\}}$ which is reducible,  and therefore is not 
isomorphic to a line.  
\smallskip 
   
 The very presence of several equivalent definitions 
 of  retracts shows how natural these objects are. 
\smallskip 

 In [7], Costa has proved that every proper retract of    
   $ K[x, y]$ (i.e., a one different  from $K$   and  $ K[x, y]$) 
has the form  $K[p]$  for some polynomial $p \in  K[x, y]$,  ~i.e., 
 is isomorphic to a polynomial $ K$-algebra in one variable. 
 A natural problem now is to  characterize somehow those polynomials $p \in  K[x, y]$ that 
 generate a retract of  $ K[x, y]$. Since the image of  a retract 
 under  any automorphism of $ K[x, y]$  is again  a retract, it would
be reasonable to characterize retracts up to  an automorphism of  $
K[x, y]$, i.e., up to a ``change of  coordinates''. 
 We give an answer to this problem  in the following 
\medskip

\noindent {\bf Theorem 3.1.}[19] Let $K[p]$ be a retract of  $ K[x,
y]$. There is 
 an automorphism $\psi$ of $ K[x, y]$ that takes the  polynomial 
$~p~$ to $~x + y \cdot q$ ~for some 
 polynomial $q = q(x,y)$. A retraction for $K[\psi(p)]$ is  given
then by $~x \to ~x + y \cdot q;   ~y \to 0$. 
\medskip 

 Geometrically, Theorem 3.1 says that   (in case $~K = 
{\Bbb C}$) every polynomial 
retraction of a plane  is a  ``parallel"  projection (sliding) on a
fiber of a coordinate polynomial (which is isomorphic to a line) 
along the fibers of another polynomial (which generates a retract of
$ K[x, y]$).  
\medskip 

 Our proof of this result is based on the   Abhyankar-Moh theorem
(see Theorem 2.2). 
\smallskip 

 Theorem 3.1 yields  another characterization of retracts of       
$ K[x,y]$:   
\medskip

\noindent {\bf Corollary 3.2.}[19] A  polynomial $~p \in  K[x, y]$ 
generates a retract of  $ K[x, y]$ 
  if and only if there is  a  polynomial mapping of $K[x,y]$ that 
takes  $~p~$ to $~x$.  The ``if" part  
is actually valid for a  polynomial algebra  in arbitrarily many
 variables. 
\smallskip 

 We also note that if  a mapping described in Corollary 3.2 
is injective, then 
$~p~$ is a  coordinate polynomial -- this follows from the
Embedding theorem of  Abhyankar and Moh [1]. 

\medskip

 Theorem 3.1  has several interesting applications, in particular, 
  to the notorious 
\medskip

\noindent {\bf Jacobian conjecture.} If for a pair of
polynomials   $p,q \in  K[x, y]$, the 
 corresponding Jacobian matrix is invertible, 
 then  $ K[p, q] =  K[x, y]$. 

\medskip

 This problem was introduced in [13], and is still unsettled. 
 For a survey and background, the reader is referred to [3].  
\smallskip

 Now we   establish a link between  retracts of   $ K[x,y]$  and the
Jacobian conjecture  by means of the following  
\medskip

\noindent {\bf Conjecture ``R''.} ~If for a pair of polynomials  
$~p,q \in  K[x, y]$, the 
 corresponding Jacobian matrix is invertible, 
 then  $K[p]$   is  a  retract of  $ K[x, y]$. 
\medskip 
 
 This statement    is formally  much   weaker than the Jacobian
conjecture since, instead of asking for   
 $~p$  ~to be a  coordinate polynomial, we  only ask for   
 $~p$  ~to generate a retract, and this property is much less 
 restrictive as can be seen from Theorem 3.1. However, the
point is that these conjectures are actually equivalent: 
\medskip 

\noindent {\bf   Theorem 3.3.}[19]    Conjecture ``R'' implies the 
Jacobian conjecture. 

\medskip

 Another application of retracts to the Jacobian conjecture 
(somewhat indirect though) is based on  the
``$\varphi^{\infty}$-trick'' familiar in combinatorial group theory
(see  [21]).  For a polynomial mapping  $\varphi : K[x,y] \to
K[x,y]$,
 ~denote by $\varphi^{\infty}(K[x,y]) = 
 \bigcap_{k=1}^{\infty}  \varphi^k(K[x,y])$ ~the {\it stable image} 
of $\varphi$. Then we have: 
\medskip

\noindent {\bf Theorem 3.4.}[19]      Let $\varphi$ be a  polynomial
mapping of $K[x,y]$. If the Jacobian 
  matrix of $\varphi$ is invertible,  then
either   $\varphi$  is   an automorphism, or  $~\varphi^{\infty}(K[x,y]) = K$. 
\smallskip 

    The proof [19] of Theorem 3.4 
 is based on a recent result of Formanek [11]. 
\smallskip 

 Obviously, if  $\varphi$ fixes  a  polynomial $p \in  K[x, y]$, ~then $~p \in \varphi^{\infty}(K[x,y])$.
 Therefore, we have: 
\medskip

\noindent {\bf Corollary 3.5.}[19] Suppose $\varphi$ is  a 
polynomial mapping of $K[x,y]$ 
 with invertible Jacobian  matrix. If $\varphi(p) = p$  ~for some
non-constant polynomial 
$p \in  K[x, y]$, ~then  $\varphi$  is   an automorphism. 
\medskip

  This yields the following  promising 
 re-formulation of the
Jacobian conjecture: if 
 $\varphi$ is  a  polynomial mapping of $K[x,y]$ 
 with invertible Jacobian  matrix, then for some automorphism 
$~\alpha$, the  mapping  $~\alpha  \cdot \varphi~$ fixes a
non-constant polynomial. 
\smallskip 

\heading 4. Some open problems \endheading 

 In this section, we have gathered  a few 
combinatorial problems about polynomial mappings that are motivated
by similar issues in combinatorial group theory. Two most 
 important problems however -- Conjectures  ``G" and  ``R" -- appear
earlier in the text (in Sections 2 and 3, respectively). 

 Throughout, $~P_n = K[x_1,...,x_n] ~$   is  the  
  polynomial algebra in  $~n$ variables, 
 $~n \ge 2$,  over a field $~K$ of characteristic $0$.
\medskip

\noindent
{\bf  (1)} [10] Is it true that every endomorphism of $~P_n~$
taking any coordinate polynomial to a coordinate one,  is actually an
automorphism? (It is true for $~n = 2$ -- see [10]).
\medskip

\noindent
{\bf  (2)} Is there a polynomial $~p \in P_n$ with the following 
property: whenever $~\varphi(p) = \psi(p)~$ for some
non-constant-valued  endomorphisms 
$~\varphi$, $\psi$ of $P_n$, it follows that $~\varphi$ = $\psi$?
(In other words, every non-constant-valued 
endomorphism of $~P_n~$ 
 is completely determined by its value on just 1 polynomial). 
\medskip

\noindent
{\bf  (3)} Suppose $~\varphi(p) = x_1~$ for some {\it
monomorphism }  (i.e., injective   endomorphism) $~\varphi~$  of
the   algebra $P_n$.  Is it true that $~p~$ is a coordinate
polynomial?  (It is true for $~n = 2$ -- see [19]). 
\medskip

\noindent
{\bf  (4)} Let $~p \in P_n$ be a polynomial  such 
  that  $~K[p]~$  is a retract of $~P_n$. Is it true that 
 $~\varphi(p) = x_1~$ for some endomorphism $~\varphi~$  of
the   algebra $P_n$ ?  (It is true for $~n = 2$ -- see [19]). 
\medskip

\noindent
{\bf  (5)} Is it true that for any endomorphism $~\varphi~$  of
the   algebra $P_n$, its stable image $\varphi^{\infty}(P_n)$ 
 is a retract of $~P_n$ ?  (The answer to this question might 
  depend on the  properties of the  ground field $~K$).

\medskip

\medskip

  Department of Mathematics, University of California, 
Santa Barbara, CA 93106 
 
{\it E-mail address\/}: shpil{\@}math.ucsb.edu

\enddocument